\documentclass[10pt,oneside]{article}
\usepackage[cp1250]{inputenc}
\usepackage[T1]{fontenc}
\usepackage{color}

\begin{document}
\newcommand{\qed}{\rule{1.5mm}{1.5mm}}
\newcommand{\proof}{\textit{Proof. }}
\newcommand{\ccon}{\rightarrowtail}
\newtheorem{theorem}{Theorem}[section]
\newtheorem{lemma}[theorem]{Lemma}
\newtheorem{remark}[theorem]{Remark}
\newtheorem{example}[theorem]{Example}
\newtheorem{corollary}[theorem]{Corollary}
\newtheorem{proposition}[theorem]{Proposition}
\newtheorem{claim}[theorem]{Claim}

\begin{center}
{\LARGE\textbf{Approximation by piecewise-regular maps}\vspace*{3mm}}\\
{\large Marcin Bilski and Wojciech Kucharz}
\end{center}
\begin{abstract}\noindent A real algebraic variety $W$ of dimension $m$ is said to be uniformly rational if
each of its points has a Zariski open neighborhood which is biregularly isomorphic to a Zariski open subset of
$\mathbf{R}^m.$ {Let $l$ be any nonnegative integer. We prove that every map of class $\mathcal{C}^l$ from a
compact subset of a real algebraic variety into a uniformly rational real algebraic variety can be approximated in
the $\mathcal{C}^l$ topology by piecewise-regular maps of class $\mathcal{C}^k,$ where $k$ is an arbitrary integer
satisfying $k\geq l$.}
{Next we derive consequences regarding algebraization of topological vector bundles.}\vspace*{1mm}\vspace*{1mm}
\end{abstract}
\section{Introduction}\label{introduction}
In this paper by a \textit{real algebraic variety} we mean a locally ringed space isomorphic to an algebraic
subset of $\mathbf{R}^m,$ for some $m,$ endowed with the Zariski topology and the sheaf of real-valued regular
functions (cf. \cite{BCR}, \cite{Ku2014}, \cite{Ku2017}). By an \textit{algebraic embedding} of a real algebraic
variety $V$ into $\mathbf{R}^n$ we mean a map $e:V\rightarrow\mathbf{R}^n$ such that $e(V)$ is a Zariski locally
closed subvariety of $\mathbf{R}^n$ and $e:V\rightarrow e(V)$ is a biregular isomorphism. Each real algebraic
variety is also equipped with the Euclidean topology induced by the standard metric on $\mathbf{R}.$ Unless
explicitly stated otherwise, all topological notions relating to real algebraic varieties refer to the Euclidean
topology.

The problem of algebraic approximation of continuous maps between real algebraic varieties has been considered by
several mathematicians (see \cite{BaFR}, \cite{BCR}, \cite{FeGh}, \cite{KuKu2} and the references therein). It is
well known that such maps can be approximated by continuous semialgebraic maps in the compact-open topology. This
is in general false if we want to approximate by regular maps instead of semialgebraic ones even for target
varieties so simple as spheres (cf. \cite{BoKu1987b}, \cite{BoKu1988}, \cite{BCR}, \cite{BoKu1987}). Therefore
various intermediate classes of maps {(more rigid than semialgebraic ones, but with better approximation
properties than regular ones)} have been investigated.

One of such classes is the class of continuous rational maps (see \cite{Ku2009}) which on nonsingular varieties
coincides with the class of regulous maps (also known as continuous hereditarily rational maps or
stratified-regular maps cf. \cite{FHMM}, \cite{KoNo}, \cite{KuKu3}). These maps have attracted a lot of attention
in recent years (see \cite{FHMM}, \cite{KoKuKu}, \cite{KoNo}, \cite{Ku2014}, \cite{KuKu2}, \cite{KuKu3},
\cite{MoJP}, \cite{Zie} and the references therein). It has turned out, for example, that every continuous map
between spheres can be approximated by regulous {ones} (see \cite{Ku2014}). However, not every continuous map from
an arbitrary compact {nonsingular} real algebraic variety into a sphere can be approximated by regulous {ones}
(see also \cite{Ku2014}).

Approximation of continuous maps from any compact {subsets of real algebraic varieties} into spheres has been
recently studied in \cite{Bi}. {The main result of \cite{Bi} says that every such map can} be approximated by
quasi-regulous maps which are obtained from regulous ones by changing signs {of the components} on some subsets of
their domains.

In the present paper we approximate maps from arbitrary compact subsets of real algebraic varieties into uniformly
rational {real algebraic} varieties (for definition see Section \ref{urrav} below). {Uniformly rational real
algebraic varieties} constitute a large class containing spheres, Grassmannians (especially interesting from the
point of view of the theory of vector bundles), rational nonsingular real algebraic surfaces and many others ({cf.
Section \ref{urrav}}). Enlarging the set of target varieties requires enlarging the class of approximating maps.
Namely, we work with piecewise-regular maps introduced in \cite{Ku2017} (see Section \ref{prm} below) the class of
which contains regulous and quasi-regulous maps mentioned above as proper subclasses (cf. \cite{Bi},
{Corollary~1}); approximating maps obtained in the present paper are neither regulous nor quasi-regulous so we do
not generalize here the main results of \cite{Ku2014} or \cite{Bi}. However, we do generalize Theorems~1.3, 1.5,
1.6 of
\cite{Ku2017} and their consequences, but not Theorem 1.4 of \cite{Ku2017}.\vspace*{2mm}\\
\textbf{Definition.} Let $V, W$ be (possibly singular) Zariski locally closed subvarieties of $\mathbf{R}^n,
\mathbf{R}^p$, respectively, $L$ a compact subset of $V,$ and $l\leq k$ nonnegative integers.

We say that a map $f:L\rightarrow W$ is of \textit{class} $\mathcal{C}^l$ (or a $\mathcal{C}^l$ \textit{map}) if
it is the restriction of some map $\mathbf{R}^n\rightarrow\mathbf{R}^p$ of class $\mathcal{C}^l$ (equivalently, if
it is the restriction of some $\mathcal{C}^l$ map $U\rightarrow\mathbf{R}^p,$ where $U$ is an open neighborhood of
$L$ in $\mathbf{R}^n$). We say that a map $f:L\rightarrow W$ is a \textit{$\mathcal{C}^l$ piecewise-regular map}
if it is both of class $\mathcal{C}^l$ and piecewise-regular.

We say that a $\mathcal{C}^l$ map $f:L\rightarrow W$ can be $\mathcal{C}^l$ \textit{approximated by
$\mathcal{C}^k$ piecewise-regular maps} if for every $\varepsilon>0$ and every $\mathcal{C}^l$ extension
$\phi=(\phi_1,\ldots,\phi_p):\mathbf{R}^n\rightarrow\mathbf{R}^p$ of $f,$ there exists a $\mathcal{C}^k$ map
$\psi=(\psi_1,\ldots,\psi_p):\mathbf{R}^n\rightarrow\mathbf{R}^p$ such that $\psi(L)\subset W,$ the restriction
$\psi|_{L}:L\rightarrow W$ is a piecewise-regular map, and
$$\big{|}\frac{\partial^{|\alpha|}{\psi}_i}{\partial
x_1^{\alpha_1}\cdots\partial x_n^{\alpha_n}}(x)-\frac{\partial^{|\alpha|}{\phi}_i}{\partial
x_1^{\alpha_1}\cdots\partial x_n^{\alpha_n}}(x)\big{|}<\varepsilon$$ for all $x\in L,$ $1\leq i\leq p,$ and
$\alpha=(\alpha_1,\ldots,\alpha_n)\in\mathbf{N}^n$ with $|\alpha|:=\alpha_1+\ldots+\alpha_n\leq l.$\vspace*{3mm}

In Section \ref{bascon} below we show that the $\mathcal{C}^l$ approximation is compatible with some topology,
which we call the $\mathcal{C}^l$ \textit{topology}, on the set $\mathcal{C}^l(L,W)$ of all $\mathcal{C}^l$~maps
from $L$ to $W.$ Thus a $\mathcal{C}^l$ map $f:L\rightarrow W$ can be $\mathcal{C}^l$ approximated by
$\mathcal{C}^k$~piecewise-regular maps if and only if every open neighborhood of $f$ in $\mathcal{C}^l(L,W)$
contains a $\mathcal{C}^k$ piecewise-regular map (see Claim \ref{claimequivalent}).

If $V, W$ are real algebraic varieties, then the notions introduced in the definition above as well as the
topology on the set $\mathcal{C}^l(L,W)$ can be defined by means of any algebraic embeddings of $V, W$ in some
$\mathbf{R}^n, \mathbf{R}^p,$ respectively, independently of the choice of the embeddings (cf. Section
\ref{ivariet}).

We point out that the $\mathcal{C}^0$ topology coincides with the usual compact-open topology. Additionally, if
$L$ is a compact $\mathcal{C}^{\infty}$ submanifold of $\mathbf{R}^n$ and the variety $W$ is nonsingular, then
$\mathcal{C}^l(L,W)$ is the space of all $\mathcal{C}^l$ maps in the sense of differential manifolds, equipped
with the compact-open $\mathcal{C}^l$ topology discussed in \cite{Hi}, p. 34 (which, by compactness of $L,$ is the
same as the Whitney topology on $\mathcal{C}^l(L,W)$).

Our main result is the following

{\begin{theorem}\label{main}Let $V, W$ be real algebraic varieties, $L$ a compact subset of $V$, and $l\leq k$
nonnegative integers. Assume that the variety $W$ is uniform\-ly rational. Then every $\mathcal{C}^l$ map
$f:L\rightarrow W$ can be $\mathcal{C}^l$ approximated by $\mathcal{C}^k$ piecewise-regular maps.
\end{theorem}}
A natural question is whether Theorem \ref{main} remains true if $V, W$ are algebraic subvarieties of some
$\mathbf{R}^n, \mathbf{R}^p,$ respectively, and $W$ is uniformly rational, but "piecewise-regular" is replaced by
"piecewise-polynomial". A polynomial map from a subset of $\mathbf{R}^n$ to a subset of $\mathbf{R}^p$ is the
restriction of a polynomial map from $\mathbf{R}^n$ to $\mathbf{R}^p.$ We obtain piecewise-polynomial maps by
substituting "polynomial" for "regular" in the definition of piecewise-regular maps. The answer to the question is
negative. Namely, by the Wood theorem (see \cite{Wood}), there are integers $q, r$ such that every polynomial map
from the unit $q$-sphere $\mathbf{S}^q$ into the unit\linebreak $r$-sphe\-re $\mathbf{S}^r$ is constant. Then,
clearly, every piecewise-polynomial map from $\mathbf{S}^q$ into $\mathbf{S}^r$ is also constant.

We point out that Theorem \ref{main} does not hold if $W$ is an arbitrary (nonrational) nonsingular real algebraic
variety (see \cite{Ku2017}, Example 1.9). Also, in general, the approximation is not possible if $k=\infty$ (see
\cite{Ku2017}, Example 1.7).

Our second result, whose presentation is postponed until Section \ref{prvb}, is an application of Theorem
\ref{main} in the theory of vector bundles. Various aspects of the problem of algebraization of topological vector
bundles on a given real algebraic variety $V$ have been studied by several mathematicians for at least the last
fifty years (see the references in \cite{BCR}, \cite{Ku2017}, \cite{KuKu2}). It has turned out that only for
exceptional varieties $V,$ every topological vector bundle on $V$ is isomorphic to an algebraic vector bundle; for
example, this is the case if $V$ is the unit $n$-sphere $\mathbf{S}^n$ (see \cite{Swa}), but need not be so if $V$
is merely assumed to be nonsingular and diffeomorphic to $\mathbf{S}^n$ (see \cite{BBK}). According to
\cite{KuKu3}, for a large class of varieties $V$, which includes all varieties homeomorphic to $\mathbf{S}^n$ and
which is known not to include all compact nonsingular varieties, every topological vector bundle on $V$ is
isomorphic to a stratified-algebraic vector bundle (or equivalently, in view of \cite{KuZie}, to a regulous vector
bundle). By Theorem 5.10 of \cite{Ku2017} every topological vector bundle on an arbitrary compact subset of any
real algebraic variety is isomorphic to a piecewise-algebraic vector bundle. In the present article we show that
the latter can be chosen of class $\mathcal{C}^k$ for any nonnegative integer $k$ (see Theorem \ref{mbund}). We do
not obtain an analogous sharpening of Theorem 5.11 of \cite{Ku2017}.

The organization of this paper is as follows. In Section \ref{DWA}, we gather preliminary material {on uniformly
rational varieties, piecewise-regular maps, $\mathcal{C}^l$~maps and the $\mathcal{C}^l$~topology.} In Section
\ref{TRZY}, the proof of Theorem \ref{main} is given. As already indicated, Section~\ref{prvb} is concerned with
vector bundles.

\section{Preliminaries}\label{DWA}
\subsection{Uniformly rational real algebraic varieties}\label{urrav}

\textbf{Definition.} Let $W$ be a real algebraic variety of dimension $n.$ A Zariski open
subset $W_0\subset W$ is said to be \textit{special} if it is biregularly isomorphic to a
Zariski open subset of $\mathbf{R}^n.$ The variety $W$ is said to be \textit{uniformly
rational} if each point of it has a special Zariski open neighborhood.\vspace*{2mm}\\
\textit{Remark.} Clearly, any uniformly rational real algebraic variety is nonsingular of pure dimension. The
question whether every nonsingular rational variety is uniformly rational remains open, see \cite{BoBo} and
\cite{Gro}, p. 885, for the discussion involving complex algebraic varieties.\vspace*{2mm}

{There are several important examples of real algebraic varieties which are known to be uniformly
rational:}\vspace*{2mm}\\
{(a) The $n$-dimensional unit sphere
$$\mathbf{S}^n=\{(x_1,\ldots,x_{n+1})\in\mathbf{R}^{n+1}:x_1^2+\ldots+x_{n+1}^2=1\}.$$ Note that
$\mathbf{S}^n\setminus\{(0,\ldots,0,1)\}$ is biregularly isomorphic to $\mathbf{R}^n$ (see \cite{BCR}, p. 76),
hence $\mathbf{S}^n$ with any point removed is isomorphic to $\mathbf{R}^n.$}\vspace*{2mm}\\
{(b) The Grassmann variety $\mathbf{G}_{k}(\mathbf{R}^n)$ of all vector subspaces of dimension $k$ of
$\mathbf{R}^n.$ Note that $\mathbf{G}_{k}(\mathbf{R}^n)$ is covered by a finite number of Zariski open sets each
of which is biregularly isomorphic to $\mathbf{R}^{(n-k)k}$ (see \cite{BCR}, p. 71 for constructing the morphisms;
analogous constructions show that for $\mathbf{F}$ equal to the field $\mathbf{C}$ of complex numbers or the field
$\mathbf{H}$ of quaternions, the variety $\mathbf{G}_k(\mathbf{F}^n)$ is a uniformly rational real algebraic
variety). }\vspace*{2mm}\\
{(c) Rational nonsingular real algebraic surfaces. This follows in principle by the Comessatti theorem (for which
see \cite{Com}, p. 257 or \cite{Kol01}, p. 206, Theorem 30 or \cite{Sil89}, Proposition 4.3). In particular, any
rational nonsingular real algebraic surface is covered by a finite number of Zariski open
subsets, each isomorphic to $\mathbf{R}^2$ (cf. \cite{Man}, Corollary 12).}\vspace*{2mm}\\
{(d) Several interesting examples can be obtained by applying the theorem saying that after blowing-ups uniformly
rational varieties remain uniformly rational (see \cite{BoBo}, \cite{Gro} for the proof in the complex setting
which also works over the field of real numbers).}
\subsection{Piecewise-regular maps}\label{prm}
Let us recall a generalization of the notion of regular map introduced in \cite{Ku2017}.\vspace*{2mm}\\
\textbf{Definition.} Let $V, W$ be real algebraic varieties, $X\subset V$ some (nonempty) subset, and $Z$ the
Zariski closure of $X$ in $V.$ A map $f:X\rightarrow W$ is said to be \textit{regular} if there is a Zariski open
neighborhood $Z_0\subseteq Z$ of $X$ and a regular map $\tilde{f}: Z_0\rightarrow W$ such that $\tilde{f}|_X=f.$\\

A \textit{stratification} of a real algebraic variety $V$ is, by definition, a finite collection of pairwise
disjoint Zariski locally closed subvarieties whose union equals~$V.$\vspace*{2mm}\\
\textbf{Definition.} Let $V, W$ be real algebraic varieties, $f:X\rightarrow W$ a continuous map defined on some
subset $X\subset V,$ and $\mathcal{S}$ a stratification of $V.$ The map $f$ is said to be \textit{piecewise
$\mathcal{S}$-regular} if for every stratum $S\in\mathcal{S}$ the restriction of $f$ to each connected component
of $X\cap S$ is a regular map (when $X\cap S$ is nonempty). Moreover, $f$ is said to be \textit{piecewise-regular}
if it is piecewise $\mathcal{T}$-regular for some stratification $\mathcal{T}$ of $V.$\vspace*{2mm}

The notion of piecewise-regular map does not depend on the ambient variety $V.$ To be precise, suppose that $V$ is
a Zariski locally closed subvariety of a real algebraic variety $V'.$ Then piecewise-regularity of $f:X\rightarrow
W$ does not depend on whether $X$ is viewed as a subset of $V$ or as that of $V'$ (see \cite{Ku2017}, p. 1546).

The following remark is an immediate consequence of the definition.\vspace*{2mm}
\\
\textit{Remark.} Let $V, W, \tilde{V}, \tilde{W}$ be real algebraic varieties and let $X\subset V$ be any subset.
Then the family of all piecewise-regular real-valued functions on $X$ constitutes a ring. Moreover, if
$f:X\rightarrow W$ is a piecewise-regular map and $\phi:\tilde{V}\rightarrow{V},$ $\psi:W\rightarrow\tilde{W}$ are
regular maps, then the composite $\psi\circ f\circ \phi|_{\phi^{-1}(X)}$ is a piecewise-regular map.\vspace*{2mm}

Let us recall the notion of \textit{nonsingular algebraic arc} (cf. \cite{KoKuKu}, \cite{Ku2017}). A subset $A$ of
a real algebraic variety $V$ is said to be a nonsingular algebraic arc if its Zariski closure $C$ in $V$ is an
algebraic curve (that is, $\mathrm{dim}(C)=1$), $A\subset C\setminus\mathrm{Sing}(C),$ and $A$ is homeomorhpic to
$\mathbf{R}$.

The following result coming from \cite{Ku2017} (Theorem 2.9) will be useful in the sequel.
\begin{theorem}\label{piecechar} Let $V, W$ be real algebraic varieties, $X\subset V$ a
semialgebraic subset,
and $f:X\rightarrow W$ a continuous semialgebraic map. Then the following conditions are equivalent:\vspace*{2mm}\\
(a) The map $f$ is piecewise-regular.\\
(b) For every nonsingular algebraic arc $A$ in $V$ with $A\subset X,$ there exists a nonempty
open subset $A_0\subset A$ such that $f|_{A_{0}}$ is a regular map.
\end{theorem}

\begin{corollary}\label{quasipiece}Let $M\subset\mathbf{R}^n$ be any semialgebraic
subset and let {$f:M\rightarrow\mathbf{R}$} be a piecewise-regular function. Let
$g:M\rightarrow\mathbf{R}$ be a continuous function such that $|f|=|g|.$ Then $g$ is a
piecewise regular function. In particular, the absolute value of every piecewise-regular
function on $M$ is a piecewise-regular function.
\end{corollary}
\proof Let $A$ be any nonsingular algebraic arc in $\mathbf{R}^n$ with $A\subset M.$ In view of Theorem
\ref{piecechar} it is sufficient to check that there exists a nonempty open subset $A_0\subset A$ such that
$g|_{A_{0}}$ is a regular function. If $A\subset f^{-1}(0),$ then $g$ is constant so it is regular. If $A$ is not
contained in the zero-set of $f$ then there is a non-empty open subset $B$ of $A$ such that $g=f$ on $B$ or $g=-f$
on $B.$ Since $f, -f$ are piecewise-regular, there is a nonempty open subset $A_0$ of $B$ such that $g$ is regular
on $A_0.$\qed
\subsection{$\mathcal{C}^l$ maps and the $\mathcal{C}^l$ topology}\label{Cltop}
In this section we provide a detailed discussion of the notions of $\mathcal{C}^l$ map,
$\mathcal{C}^l$~approximation and $\mathcal{C}^l$ topology that appear in Section \ref{introduction}. For related
constructions of topologies on map spaces we refer the reader to \cite{BaFR} and the references therein.
Henceforth $l$ stands for a nonnegative integer.
\subsubsection{Basic constructions}\label{bascon}
{Let $L$ be a nonempty compact subset of $\mathbf{R}^n$ and let $Y$ be an arbitrary subset of $\mathbf{R}^p$. A
map $f:L\rightarrow Y$ is said to be of \textit{class} $\mathcal{C}^l$ (or a $\mathcal{C}^l$ \textit{map}) if it
is the restriction of a $\mathcal{C}^l$ map $\mathbf{R}^n\rightarrow\mathbf{R}^p.$ Denote by $\mathcal{C}^l(L,Y)$
the set of all $\mathcal{C}^l$ maps from $L$ to $Y.$ Clearly,
$\mathcal{C}^{l}(L,Y)\subset\mathcal{C}^l(L,\mathbf{R}^p).$

For an arbitrary open neighborhood $U\subset\mathbf{R}^n$ of $L,$ the restriction map
$$\varrho_U:\mathcal{C}^l(U,\mathbf{R}^p)\rightarrow\mathcal{C}^l(L,\mathbf{R}^p),\mbox{ }\phi\mapsto\phi|_L$$
is surjective; here $\mathcal{C}^l(U,\mathbf{R}^p)$ is the $\mathbf{R}$-vector space of all $\mathcal{C}^l$ maps
form $U$ to $\mathbf{R}^p.$ Define the seminorm $||\cdot||:\mathcal{C}^l(U,\mathbf{R}^p)\rightarrow\mathbf{R}$ by
$$||\phi||:=\sum\sup_L\big{|}\frac{\partial^{|\alpha|} \phi_j}{\partial
x_1^{\alpha_1}\cdots\partial x_n^{\alpha_n}}\big{|}$$ for
$\phi=(\phi_1,\ldots,\phi_p)\in\mathcal{C}^l(U,\mathbf{R}^p),$ where the summation is over $j=1,\ldots,p$ and
$\alpha=(\alpha_1,\ldots,\alpha_n)\in\mathbf{N}^n$ with $|\alpha|:=\alpha_1+\ldots+\alpha_n\leq l.$ We endow
$\mathcal{C}^l(U,\mathbf{R}^p)$ with the topology induced by this seminorm. For any $\varepsilon>0,$ set
$$\mathcal{U}(\phi;U,\varepsilon):=\{\psi\in\mathcal{C}^l(U,\mathbf{R}^p):||\psi-\phi||<\varepsilon\}.$$
As $\varepsilon$ varies, the sets $\mathcal{U}(\phi;U,\varepsilon)$ constitute a base of open neighborhoods of
$\phi.$

We endow $\mathcal{C}^l(L,\mathbf{R}^p)$ with the quotient topology determined by $\varrho_U,$ which we call the
$\mathcal{C}^l$ \textit{topology}. Thus a subset $\Omega$ of $\mathcal{C}^l(L,\mathbf{R}^p)$ is open if and only
if its preimage $\varrho_U^{-1}(\Omega)$ is open in the space $\mathcal{C}^l(U,\mathbf{R}^p).$
\begin{claim}\label{claimopen}
The map $\varrho_U:\mathcal{C}^l(U,\mathbf{R}^p)\rightarrow\mathcal{C}^l(L,\mathbf{R}^p)$ is open.
\end{claim}
\proof It suffices to prove that for every $\phi\in\mathcal{C}^l(U,\mathbf{R}^p)$ and every $\varepsilon>0,$ the
set $\mathcal{U}:=\varrho_U^{-1}(\varrho_U(\mathcal{U}(\phi;U,\varepsilon)))$ is open in
$\mathcal{C}^l(U,\mathbf{R}^p).$ Given $\lambda\in\mathcal{U},$ we can choose
$\psi\in\mathcal{U}(\phi;U,\varepsilon)$ with $\varrho_U(\psi)=\varrho_U(\lambda).$ Now, setting
$\varepsilon':=||\psi-\phi||,$ we get $\mathcal{U}(\lambda;U,\varepsilon-\varepsilon')\subset\mathcal{U}.$ Indeed,
if $\mu\in\mathcal{U}(\lambda;U,\varepsilon-\varepsilon')$ and $\gamma:=\mu-\lambda+\psi,$ then
$\varrho_U(\gamma)=\varrho_U(\mu)$ and $||\gamma-\phi||\leq ||\mu-\lambda||+||\psi-\phi||<\varepsilon,$ hence
$\gamma\in\mathcal{U}(\phi;U,\varepsilon)$ and $\mu\in\mathcal{U}.$ This proves that the set $\mathcal{U}$ is
open.\qed
\begin{claim}\label{claimbasis} Let $f:L\rightarrow\mathbf{R}^p$ be a $\mathcal{C}^l$ map and let
$\phi:U\rightarrow\mathbf{R}^p$ be a $\mathcal{C}^l$ map with $\varrho_U(\phi)=f.$ Then, as $\varepsilon>0$
varies, the sets $\varrho_U(\mathcal{U}(\phi;U,\varepsilon))$ constitute a base of open neighborhoods of $f$ in
$\mathcal{C}^l(L,\mathbf{R}^p).$
\end{claim}
\proof This, in view of continuity of the map $\varrho_U,$ follows from Claim~\ref{claimopen}.\qed
\begin{claim}\label{claimindep} The $\mathcal{C}^l$ topology on $\mathcal{C}^l(L,\mathbf{R}^p)$ does not depend on
the choice of the open neighborhood $U\subset\mathbf{R}^n$ of $L.$
\end{claim}
\proof On the set $\mathcal{C}^l(L,\mathbf{R}^p)$ consider the quotient topology determined by $\varrho_U$ and the
quotient topology determined by $\varrho_{\mathbf{R}^n}.$ It suffices to show that these two topologies are
identical. To this end, choose a $\mathcal{C}^{\infty}$ function $\alpha:\mathbf{R}^n\rightarrow\mathbf{R}$ such
that $\alpha=1$ in a neighborhood of $L$ and the support of $\alpha$ is contained in $U.$ For
$\phi\in\mathcal{C}^l(U,\mathbf{R}^p)$, $\psi\in\mathcal{C}^l(\mathbf{R}^n,\mathbf{R}^p),$ and $\varepsilon>0,$ we
get
$$\varrho_{\mathbf{R}^n}(\mathcal{U}(\alpha\phi;\mathbf{R}^n,\varepsilon))=\varrho_U(\mathcal{U}(\phi;U,\varepsilon))\mbox{ and }$$
$$\varrho_{\mathbf{R}^n}(\mathcal{U}(\psi;\mathbf{R}^n,\varepsilon))=\varrho_U(\mathcal{U}(\psi|_U;U,\varepsilon)),$$
where $\alpha\phi$ is regarded as a $\mathcal{C}^l$ map on $\mathbf{R}^n.$ The proof is complete in view of
Claim~\ref{claimbasis}.\qed\vspace*{2mm}

We define the $\mathcal{C}^l$ \textit{topology} on $\mathcal{C}^l(L,Y)$ as the induced topology by regarding
$\mathcal{C}^l(L,Y)$ as the subspace of $\mathcal{C}^l(L,\mathbf{R}^p).$

For the sake of clarity, we explicitly formulate the following observation which is an immediate consequence of
Claims \ref{claimbasis}, \ref{claimindep}.
\begin{claim}\label{claimequivalent} Let $f:L\rightarrow W$ be a $\mathcal{C}^l$ map, where $W$ is a Zariski
locally closed subvariety of $\mathbf{R}^p$, and let $k\geq l$ be an integer. Then the following conditions are equivalent:\vspace*{2mm}\\
(a) The map $f$ can be $\mathcal{C}^l$ approximated by $\mathcal{C}^k$ piecewise-regular maps (in the sense
of the Definition in Section \ref{introduction} with $V=\mathbf{R}^n$).\\
(b) Every open neighborhood of $f$ in $\mathcal{C}^l(L,W)$ contains a $\mathcal{C}^k$ piecewise-regular map.\\
(c) For every $\varepsilon>0,$ there exist a $\mathcal{C}^l$ extension $\phi:U\rightarrow\mathbf{R}^p$ of $f$ and
a $\mathcal{C}^k$ map $\psi:U\rightarrow\mathbf{R}^p,$ where $U$ is an open neighborhood of $L$ in $\mathbf{R}^n,$
such that $\psi(L)\subset W,$ the restriction $\psi|_L:L\rightarrow W$ is a piecewise-regular map, and
$||\phi-\psi||<\varepsilon.$

\end{claim}
\subsubsection{Constructions involving varieties}\label{ivariet} Let $V,W$ be real algebraic varieties,
let $L$ be a compact subset of $V,$ and $k\geq l$ an integer. Fix algebraic embeddings
$e_V:V\rightarrow\mathbf{R}^n,$ $e_W:W\rightarrow\mathbf{R}^p.$ By abuse of notation, we also write
$e_V:V\rightarrow e_V(V),$ $e_W:W\rightarrow e_W(W)$ for the corresponding biregular isomorphisms. Moreover, let
$e_L:L\rightarrow e_V(L)$ denote the
restriction of $e_V.$\vspace*{2mm}\\
\textbf{Definition.} We say that a map $f:L\rightarrow W$ is of \textit{class} $\mathcal{C}^l$ (or a
$\mathcal{C}^l$ \textit{map}) if the map $e_W\circ f\circ e_L^{-1}:e_V(L)\rightarrow e_W(W)$ is of class
$\mathcal{C}^l$ in the sense of the Definition in Section~\ref{introduction}. We say that a map $f:L\rightarrow W$
is a \textit{$\mathcal{C}^l$ piecewise-regular map} if it is both of class $\mathcal{C}^l$ and piecewise-regular.

We say that a $\mathcal{C}^l$ map $f:L\rightarrow W$ can be $\mathcal{C}^l$ approximated by
$\mathcal{C}^k$~piecewise-regular maps if the map $e_W\circ f\circ e_L^{-1}:e_V(L)\rightarrow e_W(W)$ can be
$\mathcal{C}^l$ approximated by $\mathcal{C}^k$~piecewise-regular maps in the sense of the Definition in
Section~\ref{introduction}.\vspace*{2mm}

Note that a map $f:L\rightarrow W$ is piecewise-regular if and only if the map $e_W\circ f\circ
e_L^{-1}:e_V(L)\rightarrow e_W(W)$ is piecewise-regular.

The space $\mathcal{C}^l(e_V(L),e_W(W))$ of $\mathcal{C}^l$ maps endowed with the $\mathcal{C}^l$ topology is
already defined. Denoting by $\mathcal{C}^l(L,W)$ the set of all $\mathcal{C}^l$ maps from $L$ to $W,$ we get a
bijection
$$\Phi(e_V,e_W):\mathcal{C}^l(L,W)\rightarrow\mathcal{C}^l(e_V(L),e_W(W)), \mbox{ }f\mapsto e_W\circ f\circ
e_L^{-1}$$ and define the $\mathcal{C}^l$ \textit{topology} on $\mathcal{C}^l(L,W)$ so that $\Phi(e_V,e_W)$
becomes a homeomorphism.
\begin{claim}\label{claimtopindep}The set $\mathcal{C}^l(L,W)$ and the $\mathcal{C}^l$ topology on it do not depend on the choice
of the algebraic embeddings $e_V, e_W.$
\end{claim}
\proof Let $\tilde{e}_V:V\rightarrow\mathbf{R}^m, \tilde{e}_W:W\rightarrow\mathbf{R}^q$ be some algebraic
embeddings, and let $\tilde{e}_L:L\rightarrow\tilde{e}_V(L)$ denote the restriction of $\tilde{e}_V.$ Since both
maps $e_V\circ\tilde{e}_V^{-1}:\tilde{e}_V(V)\rightarrow e_V(V),$ $\tilde{e}_W\circ
e_W^{-1}:e_W(W)\rightarrow\tilde{e}_W(W)$ are biregular isomorphisms, we get a well defined homeomorphism
$$\mathcal{C}^l(e_V(L),e_W(W))\rightarrow\mathcal{C}^l(\tilde{e}_V(L),\tilde{e}_W(W)),
\mbox{ }g\mapsto \tilde{e}_W\circ e_W^{-1}\circ g\circ e_L\circ\tilde{e}_L^{-1}.$$ Clearly, this completes the
proof.\qed\vspace*{2mm}

Having defined the space $\mathcal{C}^l(L,W)$ and keeping in mind Claims \ref{claimequivalent},
\ref{claimtopindep} we see that the notion of $\mathcal{C}^l$ approximation by $\mathcal{C}^k$ piecewise-regular
maps does not depend on the choice of the algebraic embeddings $e_V, e_W,$ and Theorem \ref{main} can be restated
in the following equivalent form.
\begin{theorem}\label{maineq} Let $V, W$ be real algebraic varieties, $L$ a compact subset of $V,$ and $l\leq k$
nonnegative integers. Assume that the variety $W$ is uniformly rational. Then, for every $\mathcal{C}^l$ map
$f:L\rightarrow W$ and every  neighborhood $\Omega\subset\mathcal{C}^l(L,W)$ of $f$ in the $\mathcal{C}^l$
topology, there exists a $\mathcal{C}^k$ piecewise-regular map $g:L\rightarrow W$ that belongs to $\Omega.$
\end{theorem}

In the next two claims we point out that our definition of the space $\mathcal{C}^l(L,W)$ is compatible with the
relevant standard topological constructions.
\begin{claim}\label{claimC0} The space $\mathcal{C}^0(L,W)$ with the $\mathcal{C}^0$ topology coincides with the
space of all continuous maps from $L$ to $W$ with the compact-open topology.
\end{claim}
\proof This follows from the Tietze extension theorem for continuous functions.\qed
\begin{claim}\label{cliamCk} Assume that $e_V(L)$ is a $\mathcal{C}^{\infty}$ submanifold of $\mathbf{R}^n$
and the variety $W$ is nonsingular. Then the set $\mathcal{C}^l(L,W)$ defined in this section coincides with the
set of all $\mathcal{C}^l$ maps from $L$ to $W$ in the sense of differential manifolds. Moreover, the
$\mathcal{C}^l$ topology on $\mathcal{C}^l(L,W)$ is identical with the usual compact-open $\mathcal{C}^l$ topology
discussed in \cite{Hi}, p.34.
\end{claim}
\proof We may assume that $e_V, e_W$ are the inclusion maps. Then $L, W$ are $\mathcal{C}^{\infty}$ submanifolds
of $\mathbf{R}^n, \mathbf{R}^p,$ respectively, and using tubular neighborhoods we see that the first part of the
claim holds. Clearly, in the rest of the proof it is\linebreak\linebreak sufficient to consider the case
$W=\mathbf{R}^p$. Let $\mathcal{C}^l_{co}(L,\mathbf{R}^p)$ denote the space of all $\mathcal{C}^l$ maps from $L$
to $\mathbf{R}^p$ equipped with the compact-open $\mathcal{C}^l$ topology. Let $U$ be a tubular neighborhood of
$L$ in $\mathbf{R}^n,$ and $r:U\rightarrow L$ a $\mathcal{C}^{\infty}$ retraction.

The map $$\hat{\varrho}_U:\mathcal{C}^{l}(U,\mathbf{R}^p)\rightarrow\mathcal{C}^l_{co}(L,\mathbf{R}^p),\mbox{
}\phi\mapsto\phi|_L$$ is continuous. We will show that $\hat{\varrho}_U$ is open which, in view of Claim
\ref{claimopen}, implies that $\mathcal{C}^l(L,\mathbf{R}^p)=\mathcal{C}^l_{co}(L,\mathbf{R}^p).$ To do this, take
any open subset $\mathcal{U}$ of $\mathcal{C}^l(U,\mathbf{R}^p)$ and observe that, by Claim \ref{claimopen}, the
set $\mathcal{V}:=\hat{\varrho}_U^{-1}(\hat{\varrho}_U(\mathcal{U}))={\varrho}_U^{-1}({\varrho}_U(\mathcal{U}))$
is an open subset of $\mathcal{C}^l(U,\mathbf{R}^p).$ Moreover, we have
$\hat{\varrho}_U(\mathcal{U})=\hat{\varrho}_U(\mathcal{V}).$ Therefore, it suffices to show that
$\hat{\varrho}_U(\mathcal{V})$ is an open subset of $\mathcal{C}^l_{co}(L,\mathbf{R}^p).$

The map
$$\sigma_U:\mathcal{C}^{l}_{co}(L,\mathbf{R}^p)\rightarrow\mathcal{C}^l(U,\mathbf{R}^p),\mbox{
}\alpha\mapsto\alpha\circ r$$ is clearly continuous. Take any $f\in\hat{\varrho}_U(\mathcal{V}).$ Then
$\sigma_U(f)\in\mathcal{V}$ as $\mathcal{V}=\hat{\varrho}_U^{-1}(\hat{\varrho}_U(\mathcal{V})).$ By continuity of
$\sigma_U,$ there is an open neighborhood $\mathcal{W}$ of $f$ in $\mathcal{C}^l_{co}(L,\mathbf{R}^p)$ such that
$\sigma_U(\mathcal{W})\subset\mathcal{V}.$ Now, for every $g\in\mathcal{W},$ we have $\sigma_U(g)=g\circ
r\in\mathcal{V},$ so $\hat{\varrho}_U(g\circ r)=g\in\hat{\varrho}_U(\mathcal{V}).$ Thus, we get
$\mathcal{W}\subset\hat{\varrho}_U(\mathcal{V}),$ which proves that $\hat{\varrho}_U(\mathcal{V})$ is an open
subset of $\mathcal{C}^l_{co}(L,\mathbf{R}^p).$\qed\vspace*{2mm}

\subsubsection{A useful property of $\mathcal{C}^k$ functions}

In the proof of our main result we shall use the property that $\mathcal{C}^k$ functions which are flat on their
zero-sets remain $\mathcal{C}^k$ functions after some modifi\-cations. The class of functions defined below will
{appear in the sequel}.\vspace*{2mm}\\
\textbf{Definition.} Let $l, k$ be nonnegative integers, $l\neq 0.$ For any {open subset $U$} of $\mathbf{R}^n,$
let $\mathcal{C}^k_{l}(U)$ denote the class of all functions $v:U\rightarrow\mathbf{R}$ for which
$|v|^{\frac{1}{l}}\in\mathcal{C}^k(U).$\vspace*{2mm}

The following fact from \cite{Bi} (Lemma 3) will be useful.
\begin{lemma}\label{symmetr}Let {$U$ be an open subset} of $\mathbf{R}^n.$
Let $f\in\mathcal{C}^{k}_{l}(U)$, where $l, k$ are integers with $k\geq 1$ and $l\geq k+1.$ If
$g:U\rightarrow\mathbf{R}$ is a continuous function such that $|g(x)|=|f(x)|$ for all $x\in U,$ then
$g\in\mathcal{C}^k(U).$
\end{lemma}

\section{Proof of Theorem \ref{main}}\label{TRZY}
For any subset $A$ of $\mathbf{R}^n,$ let $\overline{A}$ denote the closure of $A$ with respect to the Euclidean
topology of $\mathbf{R}^n$ and let $\partial A$ denote the boundary of $A.$ For any two subsets $B, C$ of
$\mathbf{R}^n,$ by writing $B\subset\subset C$ we mean that $\overline{B}$ is a compact set which is contained in
$C.$

The following lemma is our main tool.
\begin{lemma}\label{partition}Let $K\subset \mathbf{R}^n$ be a compact set
and $k$ a nonnegative integer. Then for every open neighborhood $U$ of ${K}$ in $\mathbf{R}^n$ there are open
semialgebraic neighborhoods $N_1\subset\subset N_2\subset\subset U$ of ${K}$ and a piecewise-regular function
$\beta: \mathbf{R}^n\rightarrow\mathbf{R}$ of class $\mathcal{C}^k$ with $\beta(\mathbf{R}^n)\subset[0,1]$ and
with
the following properties:\vspace*{2mm}\\
\\
(1) $\partial N_1$ and $\partial N_2$ are unions of connected components of nonsingular
algebraic subvarieties of $\mathbf{R}^n$ of pure codimension $1,$\vspace*{2mm}\\
(2) $\beta|_{\mathbf{R}^n\setminus N_2}=0$ and $\beta|_{\overline{N_1}}=1.$\vspace*{2mm}\\
{In particular, all partial derivatives of $\beta$ of order from $1$ up to $k$ vanish at every point of $\partial
N_1\cup\partial N_2.$}
\end{lemma}
\proof Let $U$ be an open neighborhood of ${K}$ in $\mathbf{R}^n$. Without loss of generality we may assume that
$\overline{U}$ is compact and semialgebraic. Note that every continuous nonnegative function
$j:\overline{U}\rightarrow\mathbf{R}$ can be uniformly approximated on $\overline{U}$ by nonnegative polynomials.
Indeed, it is sufficient to approximate $\sqrt{j}$, using the Weierstrass approximation theorem, by a polynomial
$W.$ Then $W^2$ approximates~$j.$

{Applying the Weierstrass theorem and the Sard theorem do the following.} Approximate on $\overline{U}$ the
continuous function $\mathrm{dist}(\cdot, {K})$ by a nonnegative polynomial $P$ and pick
$\varepsilon_2>\varepsilon_1>0$ so that $P_1=P-\varepsilon_1, P_2=P-\varepsilon_2$ have the following properties:
$\inf_{x\in\partial U}P_1(x)>\inf_{x\in\partial U}P_2(x)>0$ and $0$ is a regular value of $P_1$ and $P_2,$ and
$P_1|_{{K}}<0,$ $P_2|_{{K}}<0.$

Now define $N_1=\{x\in U: P_1(x)<0\},$ $N_2=\{x\in U: P_2(x)<0\}$ and observe that, by the previous paragraph,
$(1)$ holds true.

Let us construct $\beta.$ {Observe that $\{P=\delta\}\cap U\neq\emptyset$ for every
$\delta\in(\varepsilon_1,\varepsilon_2).$ This is because $\sup_{x\in K}P(x)<\varepsilon_1$ and
$\inf_{x\in\partial U}P(x)>\varepsilon_2$ so, by continuity, $P$ attains all values from $(\varepsilon_1,
\varepsilon_2)$ on $U.$} {Using the Sard theorem, fix $\delta\in(\varepsilon_1,\varepsilon_2)$ such that $0$ is a
regular value of $P-\delta.$ Put $F:=(P_1\cdot P_2)^2$ and $\alpha:=\inf_{x\in\{P=\delta\}\cap U}F(x).$ Then
$\alpha>0$ as $F$ vanishes only at $x$ such that $P(x)=\varepsilon_1$ or $P(x)=\varepsilon_2$ and
$\{P=\delta\}\cap U\neq\emptyset$ is a compact set as $U$ is bounded.}

{Note that $\{F=\gamma\}\cap U\neq\emptyset$ for every $\gamma\in(0,\alpha).$ Indeed, let $\Gamma\subset U$ be an
arc connecting any point in $K$ with some point in $\partial U.$ By the previous paragraph, there are $a,b\in
\Gamma$ such that $P(a)=\varepsilon_1$ and $P(b)=\delta.$ Hence, $F(a)=0$ and $F(b)\geq\alpha$ so, by continuity,
$F$ attains all values from $(0,\alpha)$ on $\Gamma.$ Once again, using the Sard theorem, fix
$\gamma\in(0,\alpha)$ to be a regular value of $F$ and define $T=F^{-1}(\gamma).$}

Observe that for every $x\in\partial N_1$ and $y\in\partial N_2$ we have that $x,y$ are in different connected
components of $U\setminus T.$ {Indeed, suppose there is an arc $\Upsilon\subset U\setminus T$ connecting $x, y.$
We have $P(x)=\varepsilon_1$ and $P(y)=\varepsilon_2$ so there is a point $c\in\Upsilon$ such that $P(c)=\delta.$
Hence, $F(c)\geq\alpha\geq\gamma$ and $F(x)=F(y)=0.$ Consequently, there is a point $d\in\Upsilon$ with
$F(d)=\gamma,$ which means that the arc intersects $T,$ a contradiction.}

Define $G=(F-\gamma)^{2m}$ on $U$ for some nonnegative integer $m.$ Then $G^{-1}(0)=T\cap U$. Now define ${G_0}$
on $U$ by setting ${G_0}=-G$ on every connected component of $U\setminus T$ which has a nonempty intersection with
$\partial N_2$ and ${G_0}=G$ on the other connected components of $U\setminus T,$ and $G_0|_{T\cap U}=0.$ By Lemma
\ref{symmetr}, we may assume that $m$ is so large that $G_0$ is of class $\mathcal{C}^k.$ Moreover, by
Corollary~\ref{quasipiece}, in view of the fact that $G$ is piecewise-regular, we conclude that $G_0$ is also
piecewise-regular.

Next define continuous functions $$G_1^+:=-|G_0-\gamma^{2m}|+\gamma^{2m}\mbox{ and }
G_1^-:=|G_1^++\gamma^{2m}|-\gamma^{2m}.$$ Corollary \ref{quasipiece} implies that these functions are
piecewise-regular on $U$. Finally, for $i\geq 1$ define continuous functions
$$G_{i+1}^+:=-|G_i^{-}-\gamma^{2m}|+\gamma^{2m}\mbox{ and }
G_{i+1}^-:=|G_{i+1}^++\gamma^{2m}|-\gamma^{2m}.$$ Inductive application of Corollary \ref{quasipiece} proves that
these functions are piecewise-regular on $U$.

By construction, we have $(G_0)^{-1}(\kappa\cdot\gamma^{2m})\subset
(G_1^+)^{-1}(\kappa\cdot\gamma^{2m})$ and
$$(G_i^+)^{-1}(\kappa\cdot\gamma^{2m})\subset(G_i^-)^{-1}(\kappa\cdot\gamma^{2m})\mbox{ and }
(G_i^-)^{-1}(\kappa\cdot\gamma^{2m})\subset(G_{i+1}^+)^{-1}(\kappa\cdot\gamma^{2m}),$$ for
$i\geq 1$ and $\kappa\in\{-1,1\}.$ Consequently, for every $i,$ we have
${G}_i^-(x)=-\gamma^{2m}$ for $x\in\partial N_2$ and ${G}_i^{-}(x)=\gamma^{2m}$ for $x\in
\partial N_1.$

From the previous paragraph it also follows that if at some point $x\in U,$ the function $G_i^-$ is not of class
$\mathcal{C}^k,$ then $|G_i^-(x)|=\gamma^{2m}.$ Indeed, if $|G_i^-(x)|\neq\gamma^{2m},$ then, by the previous
paragraph, $|G_j^-(x)|\neq\gamma^{2m}\neq|G_j^+(x)|,$ for every $j\leq i,$ and $|G_0(x)|\neq\gamma^{2m}.$ But
$G_0$ is of class $\mathcal{C}^k$ at $x$ so $G_j^-$ is of class $\mathcal{C}^k$ at $x$ for every $j\leq i,$ by
construction.

By the fact that $\overline{U}$ is compact and again by construction, for $i$ large enough,
$|G_i^-(x)|\leq\gamma^{2m}$ for every $x\in U.$ Take such an $i$ and set $\hat{G}=G_i^-.$ From what we have just
proved we know that $\hat{G}$ is of class $\mathcal{C}^k$ on $U$ possibly outside the set
$\hat{G}^{-1}(\gamma^{2m})\cup\hat{G}^{-1}(-\gamma^{2m}).$ {We check that} for a large odd integer $l$ and a large
integer $r,$ the function
$$H=\frac{1}{(2\gamma^{2m})^{lr}}\cdot((\hat{G}-\gamma^{2m})^l+(2\gamma^{2m})^l)^r$$
is of class $\mathcal{C}^k$ on $U.$

{The fact that $H$ is of class $\mathcal{C}^k$ at every point $x_0$ with $\hat{G}(x_0)=\gamma^{2m}$ can be proved
as follows. Clearly, there is an open neighborhood $E$ of $x_0$ such that $\hat{G}-\gamma^{2m}$ is of class
$\mathcal{C}^k$ on $E\setminus\Sigma,$ where $\Sigma$ is the zero-set of $\hat{G}-\gamma^{2m}.$ Since
$\hat{G}-\gamma^{2m}$ is a continuous semialgebraic function on $E$, then, for sufficiently large $l,$ the
function $(\hat{G}-\gamma^{2m})^l$ is of class $\mathcal{C}^k$ on $E.$ To show this, it is sufficient to check
that if $l$ is large enough, then every partial derivative $\theta(x)$ of $(\hat{G}(x)-\gamma^{2m})^l$ of order
$0\leq t\leq k$ on $E\setminus\Sigma$ satisfies $\lim_{x\rightarrow x_1}{\theta(x)}=0,$ for every
$x_1\in\Sigma\cap E.$}

{Observe that $\theta(x)$ is the sum of a finite (independent of $l$) number  of terms of the form:
$(\hat{G}(x)-\gamma^{2m})^{l-j}\cdot\zeta(x)$ multiplied by a constant, where $j\leq t$ and $\zeta$ is a
continuous semialgebraic function on $E\setminus\Sigma$ independent of $l$. Then the fact that $\lim_{x\rightarrow
x_1}(\hat{G}(x)-\gamma^{2m})^{l-j}\cdot\zeta(x)=0,$ for $l$ large enough, is an immediate consequence of
Proposition 2.6.4 of \cite{BCR}.}
{Hence, the function $(\hat{G}-\gamma^{2m})^l$ is of class $\mathcal{C}^k$ and therefore $H$ is of class
$\mathcal{C}^k$ on $E.$

Similarly, for $x_0$ with $\hat{G}(x_0)=-\gamma^{2m},$ there is an open neighborhood $E$ of $x_0$ such that the
continuous semialgebraic function $(\hat{G}-\gamma^{2m})^l+(2\gamma^{2m})^{l}$ is of class $\mathcal{C}^k$ on $E$
possibly outside the zero-set of $(\hat{G}-\gamma^{2m})^l+(2\gamma^{2m})^{l}$ (recall that $l$ is odd). As before,
for $r$ large enough, $H$ is of class $\mathcal{C}^k$ on $E.$}

{Moreover, it is easy to observe that $H$ satisfies $0\leq H(x)\leq 1$ for every $x\in U$,} $H|_{\partial N_1}=1,$
$H|_{\partial N_2}=0$ and all partial derivatives of $H$ up to any prescribed order vanish at every point of
$\partial N_1\cup\partial N_2,$ {for $l, r$ large enough.}

Let us define $\beta$ by $\beta|_{N_2\setminus N_1}=H|_{N_2\setminus N_1}$ and $\beta=0$ on $\mathbf{R}^n\setminus
N_2$, and $\beta=1$ on $N_1.$ The fact that $H$ is of class $\mathcal{C}^k$ and the previous paragraph imply that
$\beta$ is of class $\mathcal{C}^k$ and satisfies (2). Clearly, $\beta$ is also semialgebraic.

It remains to check that $\beta$ is a piecewise-regular function which follows by Theorem \ref{piecechar}. Indeed,
let $A$ be a nonsingular algebraic arc in $\mathbf{R}^n.$ First assume that $A\cap
(N_2\setminus\overline{N_1})\neq\emptyset.$ Then, by construction of $H,$ there is an open nonempty subset $A_0$
of $A$ contained in $N_2\setminus\overline{N_1}$ such that $\beta|_{A_0}=H|_{A_0}$ is a regular function. If
$A\cap (N_2\setminus\overline{N_1})=\emptyset,$ then there is an open nonempty subset $A_0$ of $A$ contained in
either $\overline{N_1}$ or $\mathbf{R}^n\setminus N_2.$ Then
$\beta|_{A_0}$ is constant, hence regular. Now the claim follows immediately.\qed\vspace*{2mm}\\
\textbf{Proof of Theorem \ref{main}.} By Section \ref{ivariet}, we may assume that $V, W$ are Zariski locally
closed subvarieties of some $\mathbf{R}^n, \mathbf{R}^p,$ respectively. Clearly, it is sufficient to consider the
case $V=\mathbf{R}^n.$ Let $m$ denote the dimension of $W.$ Fix $f\in\mathcal{C}^l(L,W)$ and take some extension
of $f$ belonging to $\mathcal{C}^l(\mathbf{R}^n,\mathbf{R}^p).$ The extension will also be denoted by $f.$ Since
$W$ has a tubular neighborhood in $\mathbf{R}^p,$ we may assume that $f(\Omega)\subset W$ for some open
neighborhood $\Omega$ of $L$ in $\mathbf{R}^n.$ To complete the proof it is sufficient to show that for every
$\varepsilon>0,$ there is a $\mathcal{C}^k$ piecewise-regular map $g:T\rightarrow W,$ where $T$ is an open
neighborhood of $L$ in $\mathbf{R}^n,$ satisfying $||g-f|_T||<\varepsilon$ (cf. Claim \ref{claimequivalent} with
$U=T$).

Let $c$ be the least positive integer such that the set $f(L)$ is contained in the union of $c$ special Zariski
open subsets of $W.$ The proof is by induction on $c.$

If $c=1,$ then there is a Zariski open subset $E$ of $W$ such that $f(L)\subset E$ and there is a biregular
isomorphism $\phi:D\rightarrow E,$ where $D$ is a Zariski open subset of $\mathbf{R}^m.$ Shrinking the open
neighborhood $\Omega$ of $L$ if necessary, we get $f(\Omega)\subset E,$ and hence $f|_{\Omega}=\phi\circ h$ for
some {$\mathcal{C}^l$ map $h:\Omega\rightarrow D\subset\mathbf{R}^m.$ Now it is sufficient to approximate $h,$
using the Weierstrass approximation theorem, by a polynomial map $\tilde{h}:\mathbf{R}^n\rightarrow\mathbf{R}^m$
in such a way that the following holds: the restrictions to $L$ of all partial derivatives of order up to $l$ of
the components of $\tilde{h}$ are as close to the corresponding restrictions of the partial derivatives of the
components of $h$ as we wish. Then the map $g$ equal to $\phi\circ\tilde{h}$ in some neighborhood of $L$ in
$\mathbf{R}^n$ is a $\mathcal{C}^k$ piecewise-regular map approximating $f.$}

Let $c>1$ and let $\{E_1,\ldots,E_c\}$ be a family of special subsets of $W$ such that $f(L)\subset
E_1\cup\ldots\cup E_c.$ {Note that, there is an open bounded semialgebraic neighborhood $T$ of $L$ in
$\mathbf{R}^n$ with $f(\overline{T})\subset E_1\cup\ldots\cup E_c.$ } Then the compact set
$K:=(f|_{\overline{T}})^{-1}(W\setminus E_c)$ has an open bounded semialgebraic neighborhood $U$ in $\mathbf{R}^n$
with $f(\overline{U})\subset E_1\cup\ldots\cup E_{c-1}.$

By Lemma \ref{partition}, there are open semialgebraic neighborhoods $N_1\subset\subset N_2\subset\subset U$ of
$K$ and a piecewise-regular function $\beta: \mathbf{R}^n\rightarrow\mathbf{R}$ of class $\mathcal{C}^k$ such that
$\beta(\mathbf{R}^n)\subset[0,1]$ and
the following conditions hold:\\
\\
(1) $\partial N_1$ and $\partial N_2$ are unions of connected components of nonsingular
algebraic subvarieties of $\mathbf{R}^n$ of pure codimension $1,$\vspace*{2mm}\\
(2) $\beta|_{\mathbf{R}^n\setminus N_2}=0$ and $\beta|_{\overline{N_1}}=1.$\vspace*{2mm}

Let $A\subset\subset N_1$ be an open neighborhood of ${K}.$ Note that $\overline{N_2}\setminus
N_1\subset(\overline{U}\setminus{A})=:B$ and $f(B\cap\overline{T})\subset E_c.$

Since $f(\overline{U})\subset E_1\cup\ldots\cup E_{c-1},$ then, by the induction hypothesis, there is a
$\mathcal{C}^k$ piecewise-regular map $f_1:\overline{U}\rightarrow W$ approximating $f|_{\overline{U}}$ in the
$\mathcal{C}^l$ topology. We may assume that the approximation is close enough to ensure that
$f_1(B\cap\overline{T})\subset E_c.$ Since $E_c$ is a special subset of $W,$ there is a biregular map $\phi:
D\rightarrow E_c,$ where $D$ is a Zariski open subset of $\mathbf{R}^m.$ By the inclusion
$f_1(B\cap\overline{T})\subset E_c,$ we have a $\mathcal{C}^k$ piecewise-regular map
$h_1:B\cap\overline{T}\rightarrow D$ such that $f_1|_{B\cap\overline{T}}=\phi\circ h_1.$

By the definition of $K$ and the choice of $A,$ we have $f(\overline{T}\setminus A)\subset E_c.$ Consequently, as
above, there is a {$\mathcal{C}^l$ map $h:\overline{T}\setminus A\rightarrow D$ such that
$f|_{\overline{T}\setminus A }=\phi\circ h.$}

{Now approximate $h,$ using the Weierstrass approximation theorem, by a polynomial map
${h_2}:\mathbf{R}^n\rightarrow\mathbf{R}^m$ in such a way that the following holds: the restrictions to
$\overline{T}\setminus A$ of all partial derivatives of order up to $l$ of the components of ${h}_2$  are as close
to the corresponding restrictions of the partial derivatives of the components of $h$ as we wish. Then the map
$\phi\circ{h}_2|_{\overline{T}\setminus A}$ is a $\mathcal{C}^k$ piecewise-regular map approximating
$f|_{\overline{T}\setminus A}$ in the $\mathcal{C}^l$ topology.}

Note that $B\cap\overline{T}\subset \overline{T}\setminus A$ and observe that $h_1$ and $h_2$ are close to each
other on $B\cap\overline{T}.$ Therefore the formula $\tilde{h}=\beta\cdot h_1+(1-\beta)\cdot h_2$ gives a
$\mathcal{C}^{k}$ piecewise-regular map $\tilde{h}:B\cap\overline{T}\rightarrow D$ close to
$h_1|_{B\cap\overline{T}}$ and to $h_2|_{B\cap\overline{T}},$ and to $h|_{B\cap\overline{T}}.$

Finally, let us define a semialgebraic map ${g}=({g}_1,\ldots,{g}_p):T\rightarrow W\subset\mathbf{R}^p$ by:
$${g}|_{T\setminus{N_2}}:=\phi\circ h_2|_{T\setminus N_2}\mbox{ and }{g}|_{T\cap (N_2\setminus
\overline{N_1})}:=\phi\circ\tilde{h}|_{T\cap (N_2\setminus \overline{N_1})}\mbox{ and }{g}|_{T\cap
\overline{N_1}}:=f_1|_{T\cap\overline{N_1}}$$ and let us show that ${g}$ is a $\mathcal{C}^k$ piecewise-regular
map and that the seminorm $||g-f|_T||$ is small.

Clearly, all partial derivatives of order up to $l$ of the components of ${g}$ approximate the corresponding
partial derivatives of the components of $f$ on every set of the family $\{T\setminus\overline{N_2},
T\cap(N_2\setminus\overline{N_1}), T\cap N_1\}.$ To show that ${g}$ is a $\mathcal{C}^k$ map and that the seminorm
$||g-f|_T||$ is small it remains to check that for every $\alpha\in\mathbf{N}^n$ with
$|\alpha|:=\alpha_1+\ldots+\alpha_n\leq k$ and every $i=1,\ldots,p,$ the functions
$$\frac{\partial^{|\alpha|} {g}_i}{\partial x_1^{\alpha_1}\cdots\partial x_n^{\alpha_n}}|_{T\setminus
\overline{N_2}},\mbox{ }\frac{\partial^{|\alpha|} {g}_i}{\partial x_1^{\alpha_1}\cdots\partial
x_n^{\alpha_n}}|_{T\cap(N_2\setminus\overline{N_1})},\mbox{ } \frac{\partial^{|\alpha|} {g}_i}{\partial
x_1^{\alpha_1}\cdots\partial x_n^{\alpha_n}}|_{T\cap N_1}$$ can be glued along $T\cap(\partial N_1\cup\partial
N_2)$ to constitute a continuous function on $T.$

{To do this, observe that, by the properties of $\beta$ (cf. Lemma \ref{partition}), at every $b\in\partial
N_2\cap T$ (resp. at every $a\in\partial N_1\cap T$), the corresponding partial derivatives of the components of
$\tilde{h}$ and of $h_2$ (resp. of $h_1$) are equal up to order $k.$ Therefore, the corresponding partial
derivatives of the components of $\phi\circ\tilde{h}$ and of $\phi\circ h_2$ (resp. of $\phi\circ h_1$) can be
glued along $\partial N_2\cap T$ (resp. $\partial N_1\cap T$). Now the claim follows immediately.

To complete the proof it is sufficient to show that ${g}$ is a piecewise-regular map. Here we shall use Theorem
\ref{piecechar}. Let $A$ be a nonsingular algebraic arc in $\mathbf{R}^n$ with $A\subset T.$ If $A\cap
(T\setminus(\partial N_1\cup\partial N_2))\neq\emptyset,$ then there is an open subset $A_0$ of $A$ such that
$A_0$ is contained in $T\setminus\overline{N_2}$ or in $T\cap({N_2\setminus\overline{N_1}})$ or in $T\cap N_1$ and
then the claim follows by Theorem \ref{piecechar} and the definition of ${g}.$ If $A\cap (T\setminus(\partial
N_1\cup\partial N_2))=\emptyset,$ then there is an open subset $A_0$ of $A$ such that $A_0$ is contained in
$T\cap\partial{N_2}\subset T\setminus N_2$ or in $T\cap\partial N_1\subset T\cap\overline{N_1}$ and again the
claim is a direct consequence of Theorem \ref{piecechar} and the definition of ${g}.$ \qed
\section{$\mathcal{C}^k$ piecewise-algebraic vector bundles}\label{prvb}
Piecewise-algebraic vector bundles have been introduced in \cite{Ku2017} to which we refer the reader for details.
Before stating the main result of this section we recall some terminology and facts from \cite{BCR}, \cite{Hi} and
\cite{Ku2017}.

Let $\mathbf{F}$ denote $\mathbf{R}, \mathbf{C}$ or the field $\mathbf{H}$ of quaternions. All $\mathbf{F}$-vector
spaces are assumed to be left vector spaces. This plays a role if $\mathbf{F}=\mathbf{H}$ since the quaternions
are noncommutative. Let $\xi$ be a topological $\mathbf{F}$-vector bundle over a topological space $X$. By
$E(\xi)$ we denote the total space of $\xi$ and by $p(\xi): E(\xi)\rightarrow X$ the bundle projection. The fiber
of $\xi$ over a point $x\in X$ is $E(\xi)_x=p(\xi)^{-1}(x).$

For any nonnegative integer $m$, let $\varepsilon^m_X(\mathbf{F})$ denote the standard product $\mathbf{F}$-vector
bundle on $X$ with total space $X\times\mathbf{F}^m.$ If $\xi$ is a topological $\mathbf{F}$-vector subbundle of
$\varepsilon_X^m(\mathbf{F}),$ then $\varepsilon_X^m(\mathbf{F})=\xi\oplus\xi^{\perp},$ where $\xi^{\perp}$ is the
orthogonal complement of $\xi$ with respect to the standard inner product on $\mathbf{F}^m.$ Then the orthogonal
projection $\rho_{\xi}:\varepsilon_X^m(\mathbf{F})\rightarrow\varepsilon_X^m(\mathbf{F})$ onto $\xi$ is a
topological morphism of $\mathbf{F}$-vector bundles.

Let $V$ be a real algebraic variety. Then $V\times\mathbf{F}^m$ can also be regarded as a real algebraic variety.
By an algebraic $\mathbf{F}$-vector bundle on $V$ we mean an algebraic vector subbundle of
$\varepsilon_V^m(\mathbf{F})$ for some $m$ (cf. \cite{BCR}, Chapters 12 and 13).

Let $X$ be a subspace of a topological space $Y$ and $\psi:\theta\rightarrow\omega$ a topological morphism of
topological $\mathbf{F}$-vector bundles on $Y.$ We let $\psi_X:\theta|_X\rightarrow\omega|_X$ denote the
restriction morphism defined by $\psi_X(v)=\psi(v)$ for all $v\in E(\theta|_X).$

The following generalization of the notion of algebraic vector bundle is taken from
\cite{Ku2017}.\vspace*{2mm}\\
\textbf{Definition.} Let $V$ be a real algebraic variety, $X\subset V$ some nonempty subset and $Z$ the Zariski
closure of $X$ in $V.$

An \textit{algebraic $\mathbf{F}$-vector bundle $\xi$ on} $X$ is a topological $\mathbf{F}$-vector subbundle of
$\varepsilon_{X}^m(\mathbf{F}),$ for some $m,$ for which there exist a Zariski open neighborhood $Z_0\subset Z$ of
$X$ and an algebraic $\mathbf{F}$-vector subbundle $\tilde{\xi}$ of $\varepsilon^m_{Z_0}(\mathbf{F})$ with
$\tilde{\xi}|_{X}=\xi.$ Then $\xi$ is also said to be an \textit{algebraic $\mathbf{F}$-vector subbundle} of
$\varepsilon_X^m(\mathbf{F}).$ The pair $(Z_0,\tilde{\xi})$ is said to be an algebraic extension of $\xi.$

If $\xi$, $\eta$ are algebraic $\mathbf{F}$-vector bundles on $X,$ then an \textit{algebraic morphism}
$\phi:\xi\rightarrow\eta$ is a topological morphism such that there are algebraic extensions $(Z_0,\tilde{\xi}),$
$(Z_0,\tilde{\eta})$ of $\xi, \eta,$ respectively, and an algebraic morphism
$\tilde{\phi}:\tilde{\xi}\rightarrow\tilde{\eta}$ with $\tilde{\phi}_X=\phi.$

The following notion is also taken from \cite{Ku2017}. \vspace*{2mm}\\
\textbf{Definition.} Let $V$ be a real algebraic variety, $X\subset V$ some subset, and $\mathcal{S}$ a
stratification of $V.$

A \textit{piecewise $\mathcal{S}$-algebraic $\mathbf{F}$-vector bundle $\xi$ on $X$} is a topological
$\mathbf{F}$-vector subbundle of $\varepsilon_X^m(\mathbf{F}),$ for some $m,$ such that for every stratum
$S\in\mathcal{S}$ and each connected component $\Sigma$ of $X\cap S$ the restriction $\xi|_{\Sigma}$ is an
algebraic $\mathbf{F}$-vector subbundle of $\varepsilon_{\Sigma}^m(\mathbf{F}).$ In that case, $\xi$ is said to be
a \textit{piecewise $\mathcal{S}$-algebraic $\mathbf{F}$-vector subbundle} of $\varepsilon_X^m(\mathbf{F}).$

If $\xi, \eta$ are piecewise $\mathcal{S}$-algebraic $\mathbf{F}$-vector bundles on $X,$ then a \textit{piecewise
$\mathcal{S}$-algebraic morphism} $\phi:\xi\rightarrow\eta$ is a topological morphism such that for every stratum
$S\in\mathcal{S}$ and each connected component $\Sigma$ of $S\cap X,$ the restriction
$\phi_{\Sigma}:\xi|_{\Sigma}\rightarrow\eta|_{\Sigma}$ is an algebraic morphism.

A \textit{piecewise-algebraic $\mathbf{F}$-vector bundle on $X$} is a piecewise $\mathcal{T}$-algebraic
$\mathbf{F}$-vector bundle on $X$ for some stratification $\mathcal{T}$ of $V.$

If $\xi$ and $\eta$ are piecewise-algebraic $\mathbf{F}$-vector bundles on $X,$ then a \textit{piecewise-algebraic
morphism} $\phi:\xi\rightarrow\eta$ is a piecewise $\mathcal{T}$-algebraic morphism for some stratification
$\mathcal{T}$ of $V$ such that both $\xi$ and $\eta$ are piecewise $\mathcal{T}$-algebraic $\mathbf{F}$-vector
bundles on $X.$\vspace*{2.5mm}

In what follows $k$ denotes a nonnegative integer. The notion of a bundle of class $\mathcal{C}^k$ on a smooth
manifold has been discussed in \cite{Hi}. We generalize it to the case where the base space is a compact subset
$L$ of $\mathbf{R}^n$ with the Euclidean
topology induced from $\mathbf{R}^n.$\vspace*{2mm}\\
\textbf{Definition.} A topological $\mathbf{F}$-vector bundle $\xi$ on $L$ is said to be of \textit{class}
$\mathcal{C}^k$ if $\xi$ is a subbundle of $\varepsilon_L^m(\mathbf{F}),$ for some $m,$ such that there exist an
open neighborhood $U$ of $L$ in $\mathbf{R}^n$ and an $\mathbf{F}$-vector subbundle $\hat{\xi}$ of
$\varepsilon_U^m(\mathbf{F})$ of class $\mathcal{C}^k$ satisfying $\xi=\hat{\xi}|_L.$ In that case, $\xi$ is
called a subbundle of $\varepsilon_L^m(\mathbf{F})$ of class $\mathcal{C}^k.$ The pair $(U,\hat{\xi})$ is called a
$\mathcal{C}^k$ \textit{extension of} $\xi.$ In particular, $\varepsilon^m_L(\mathbf{F})$ is an
$\mathbf{F}$-vector bundle of class $\mathcal{C}^k$ on $L.$

If $\xi, \eta$ are $\mathbf{F}$-vector bundles of class $\mathcal{C}^k$ on $L,$ then a $\mathcal{C}^k$
\textit{morphism} $\phi:\xi\rightarrow\eta$ is a topological morphism such that there exist $\mathcal{C}^k$
extensions $(U,\hat{\xi}), (U,\hat{\eta})$ of $\xi, \eta,$ respectively, and a $\mathcal{C}^k$ morphism
$\hat{\phi}:\hat{\xi}\rightarrow\hat{\eta}$ with $\hat{\phi}|_L=\phi.$\vspace*{2.5mm}

The notions introduced in the definition above can be extended in a natural way to the setting involving real
algebraic varieties. For the sake of clarity, we first recall some notation.

Let $f:X\rightarrow Y$ be a continuous map of topological spaces. If $\theta$ is a topological $\mathbf{F}$-vector
subbundle of $\varepsilon^m_Y(\mathbf{F}),$ then the pullback $f^{*}\theta$ is a subbundle of
$\varepsilon_X^m(\mathbf{F})$ with total space
$$E(f^{*}\theta)=\{(x,v)\in X\times\mathbf{F}^m:(f(x),v)\in E(\theta)\}.$$
If $\phi:\theta\rightarrow\omega$ is a morphism, where $\omega$ is an $\mathbf{F}$-vector subbundle of
$\varepsilon_Y^p(\mathbf{F}),$ then the pullback morphism $f^{*}\phi:f^{*}\theta\rightarrow f^{*}\omega$ is
defined by $(f^{*}\phi)(x,v)=(f(x), \pi(\phi(x,v))),$ where $(x,v)\in E(f^{*}\xi)$ and
$\pi:X\times\mathbf{F}^p\rightarrow\mathbf{F}^p$ is the canonical projection.\\
\textbf{Definition.} Let $V$ be a real algebraic variety and let $L$ be a compact subset of $V.$ Fix an algebraic
embedding $e:V\rightarrow\mathbf{R}^n$ and denote by $e_L:L\rightarrow e(L)$ the restriction of $e.$

A topological $\mathbf{F}$-vector bundle $\xi$ on $L$ is said to be of \textit{class} $\mathcal{C}^k$ if $\xi$ is
a subbundle of $\varepsilon_L^m(\mathbf{F})$, for some $m,$ such that the pullback $\mathbf{F}$-vector bundle
$(e_L^{-1})^{*}\xi$ is a subbundle of $\varepsilon^m_{e(L)}(\mathbf{F})$ of class $\mathcal{C}^k.$

If $\xi, \eta$ are $\mathbf{F}$-vector bundles of class $\mathcal{C}^k$ on $L,$ then a $\mathcal{C}^k$
\textit{morphism} $\phi:\xi\rightarrow\eta$ is a topological morphism such that the pullback morphism
$(e_L^{-1})^{*}\phi:(e_L^{-1})^{*}\xi\rightarrow(e_L^{-1})^{*}\eta$ is of class $\mathcal{C}^k.$\vspace*{3mm}

The notions just defined do not depend on the choice of the algebraic embedding $e:V\rightarrow
e(V)\subset\mathbf{R}^n.$ This assertion readily follows since if
$\tilde{e}:V\rightarrow\tilde{e}(V)\subset\mathbf{R}^q$ is another algebraic embedding, then the map
$e\circ\tilde{e}^{-1}:\tilde{e}(V)\rightarrow e(V)$ is a biregular isomorphism.

It is clear that if $\xi$ is an $\mathbf{F}$-vector subbundle of $\varepsilon_L^m(\mathbf{F})$ of class
$\mathcal{C}^k,$ then the orthogonal projection
$\rho_{\xi}:\varepsilon_L^m(\mathbf{F})\rightarrow\varepsilon_L^m(\mathbf{F})$ onto $\xi$ is a
$\mathcal{C}^k$ morphism of $\mathbf{F}$-vector bundles of class $\mathcal{C}^k$.

Note that the category of $\mathbf{F}$-vector bundles of class $\mathcal{C}^0$ on $L$ defined in this section is
equivalent to the category of topological $\mathbf{F}$-vector bundles on $L$ (see \cite{HuD}, p. 31, Proposition 5.8
and \cite{Hi}, p. 92, Exercise 1).\vspace*{2mm}\\
\textbf{Definition.} Let $L$ be a compact subset of a real algebraic variety $V$ and $\mathcal{S}$ a
stratification of $V.$ A \textit{$\mathcal{C}^k$ piecewise $\mathcal{S}$-algebraic $\mathbf{F}$-vector bundle on}
$L$ is a piecewise $\mathcal{S}$-algebraic $\mathbf{F}$-vector subbundle of $\varepsilon_L^m(\mathbf{F}),$ for
some $m,$ which is of class $\mathcal{C}^k$.

If $\xi, \eta$ are $\mathcal{C}^k$ piecewise $\mathcal{S}$-algebraic $\mathbf{F}$-vector bundles on $L$ then a
\textit{$\mathcal{C}^k$ piecewise $\mathcal{S}$-algebraic morphism} $\phi:\xi\rightarrow\eta$ is a morphism both
in the category of $\mathbf{F}$-vector bundles of class $\mathcal{C}^k$ and the category of piecewise
$\mathcal{S}$-algebraic $\mathbf{F}$-vector bundles.

A \textit{$\mathcal{C}^k$ piecewise-algebraic $\mathbf{F}$-vector bundle on} $L$ is a $\mathcal{C}^k$ piecewise
$\mathcal{T}$-algebraic $\mathbf{F}$-vector bundle on $L$ for some stratification $\mathcal{T}$ of $V$.

If $\xi, \eta$ are $\mathcal{C}^k$ piecewise-algebraic $\mathbf{F}$-vector bundles on $L$ then a
\textit{$\mathcal{C}^k$ piecewise-algebraic morphism} $\phi:\xi\rightarrow\eta$ is a $\mathcal{C}^k$ piecewise
$\mathcal{T}$-algebraic morphism, for some stratification $\mathcal{T}$ of $V$ such that both $\xi$ and $\eta$ are
$\mathcal{C}^k$ piecewise $\mathcal{T}$-algebraic $\mathbf{F}$-vector bundles on $L.$ \vspace*{2.5mm}

Define $\gamma_r(\mathbf{F}^m)$ to be the tautological $\mathbf{F}$-vector bundle on $\mathbf{G}_r(\mathbf{F}^m).$
The bundle $\gamma_r(\mathbf{F}^m)$ is an algebraic $\mathbf{F}$-vector subbundle of
$\varepsilon^m_{\mathbf{G}_r(\mathbf{F}^m)}(\mathbf{F}).$ Let $\mathbf{G}(\mathbf{F}^m)$ be the disjoint union of
the $\mathbf{G}_r(\mathbf{F}^m),$ $0\leq r\leq m.$ We denote by $\gamma(\mathbf{F}^m)$ the algebraic
$\mathbf{F}$-vector subbundle of $\varepsilon^m_{\mathbf{G}(\mathbf{F}^m)}(\mathbf{F})$ whose restriction to
$\mathbf{G}_r(\mathbf{F}^m)$ is $\gamma_r(\mathbf{F}^m)$, for $0\leq r\leq m.$

Let $Y$ be a topological space and $\xi$ a topological $\mathbf{F}$-vector subbundle of
$\varepsilon_Y^m(\mathbf{F}).$ Then the map $f_{\xi}:Y\rightarrow\mathbf{G}(\mathbf{F}^m)$ defined by
$E(\xi)_y=\{y\}\times f_{\xi}(y),$ for all $y\in Y,$ is continuous and $\xi=f_{\xi}^{*}\gamma(\mathbf{F}^m).$ We
call $f_{\xi}$ the \textit{classifying map} for $\xi.$

It follows immediately from the definition and \cite{Hi} that a topological $\mathbf{F}$-vector subbundle of
$\varepsilon_L^m(F),$ where $L$ is a compact subset of a real algebraic variety, is of class $\mathcal{C}^k$ if
and only if the classifying map $f_{\xi}:L\rightarrow\mathbf{G}(\mathbf{F}^m)$ is a $\mathcal{C}^k$ map. Combining
this fact with Proposition 5.6 of \cite{Ku2017} we obtain
\begin{proposition}\label{bundle1}Let $V$ be a real algebraic variety, $L\subset V$ a compact subset, $\mathcal{S}$
a stratification of $V,$ and $\xi$ a topological $\mathbf{F}$-vector subbundle of $\varepsilon_L^m(\mathbf{F})$
for some nonnegative integer $m.$ Then the following conditions are
equivalent\vspace*{2mm}:\\
(a) $\xi$ is a $\mathcal{C}^k$ piecewise $\mathcal{S}$-algebraic $\mathbf{F}$-vector subbundle of
$\varepsilon_L^m(\mathbf{F}).$\\
(b) The classifying map $f_{\xi}:L\rightarrow\mathbf{G}(\mathbf{F}^m)$ for $\xi$ is a $\mathcal{C}^k$ piecewise
$\mathcal{S}$-regular map.
\end{proposition}

The following proposition is a variant of Proposition 5.8 of \cite{Ku2017} in the category of vector bundles of
class $\mathcal{C}^k$.
\begin{proposition}\label{bundle2} Let $V$ be a real algebraic variety, $L\subset V$ a compact subset, and $\mathcal{S}$ a
stratification of $V$. Let $\xi, \eta$ be $\mathcal{C}^k$ piecewise $\mathcal{S}$-algebraic $\mathbf{F}$-vector
bundles on $L$ that are topologically isomorphic. Then $\xi$ and $\eta$ are also isomorphic in the category of
$\mathcal{C}^k$ piecewise $\mathcal{S}$-algebraic $\mathbf{F}$-vector bundles on $L.$
\end{proposition}
\proof We follow the proof of Proposition 5.8 of \cite{Ku2017}. The bundle $\xi$ (resp. $\eta$) is a piecewise
$\mathcal{S}$-algebraic $\mathbf{F}$-vector subbundle of $\varepsilon_L^p(\mathbf{F})$ (resp.
$\varepsilon_L^q(\mathbf{F})$) of class $\mathcal{C}^k$ for some $p$ (resp. $q$). Since
$\varepsilon_L^p(\mathbf{F})=\xi\oplus\xi^{\perp}$ and $\varepsilon_L^q(\mathbf{F})=\eta\oplus\eta^{\perp},$ there
exists a topological morphism $\phi:\varepsilon_L^p(\mathbf{F})\rightarrow\varepsilon_L^q(\mathbf{F})$ which
transforms $\xi$ onto $\eta.$ Let $A=A_{\phi}:L\rightarrow {\mathrm{Mat}}_{q,p}(\mathbf{F})$ be the matrix
representation of $\phi$ (cf. \cite{Ku2017}, Section 4.2). By the Weierstrass approximation theorem there is a
regular map $B:L\rightarrow \mathrm{Mat}_{q,p}(\mathbf{F})$ close to $A.$ Then $\psi:
\varepsilon_L^p(\mathbf{F})\rightarrow\varepsilon_L^q(\mathbf{F}),$ defined by
$$\psi(x,v)=(x, B(x)(v))\mbox{ for }(x,v)\in L\times\mathbf{F}^p,$$
is an algebraic morphism.

By the fact that $\eta$ is a $\mathcal{C}^k$ piecewise $\mathcal{S}$-algebraic $\mathbf{F}$-vector bundle on $L$
and by Lemma 5.3 of \cite{Ku2017}, the orthogonal projection
$\rho_{\eta}:\varepsilon_L^q(\mathbf{F})\rightarrow\varepsilon_L^q(\mathbf{F})$ onto $\eta$ is a $\mathcal{C}^k$
piecewise $\mathcal{S}$-algebraic morphism. Therefore,
$\rho_{\eta}\circ\psi:\varepsilon_L^p(\mathbf{F})\rightarrow\varepsilon_L^q(\mathbf{F})$ is a $\mathcal{C}^k$
piecewise $\mathcal{S}$-algebraic morphism which transforms $\xi$ onto $\eta.$ Consequently, the morphism
$\sigma:\xi\rightarrow\eta$ determined by $\rho_{\eta}\circ\psi$ is bijective and $\mathcal{C}^k$ piecewise
$\mathcal{S}$-algebraic and its inverse is of class $\mathcal{C}^k$. By Lemma 5.2 of \cite{Ku2017} we conclude
that $\sigma$ is a $\mathcal{C}^k$ piecewise $\mathcal{S}$-algebraic isomorphism.\qed\vspace*{0mm}\\

The following consequence of Theorem \ref{main} and Propositions \ref{bundle1}, \ref{bundle2} is the main result
of this section.
\begin{theorem}\label{mbund} Let $V$ be a real algebraic variety, $L\subset V$ a compact subset, and $l\leq k$
nonnegative integers. Then each $\mathbf{F}$-vector bundle on $L$ of class $\mathcal{C}^l$ is $\mathcal{C}^l$
isomorphic to a $\mathcal{C}^k$ piecewise-algebraic $\mathbf{F}$-vector bundle on $L.$ The latter bundle is
uniquely determined up to $\mathcal{C}^k$ piecewise-algebraic isomorphism.
\end{theorem}
\proof Let $\xi$ be an $\mathbf{F}$-vector bundle on $L$ of class $\mathcal{C}^l.$ Then there are a positive
integer $m$ and a continuous map $f:L\rightarrow\mathbf{G}(\mathbf{F}^m)$ such that $\xi$ is topologically
isomorphic to the pullback $f^{*}\gamma(\mathbf{F}^m)$ (cf. \cite{HuD}, Chapter 3, Proposition 5.8). Without loss
of generality, we may assume that $f(L)\subset\mathbf{G}_r(\mathbf{F}^m)$ for some $0\leq r\leq m.$

Recall that $\mathbf{G}_r(\mathbf{F}^m)$ is a uniformly rational real algebraic variety. Then, by Theorem
\ref{main} with $W=\mathbf{G}_r(\mathbf{F}^m),$ we obtain that $f$ is homotopic to a
$\mathcal{C}^k$~piecewise-regular map $g:L\rightarrow\mathbf{G}_r(\mathbf{F}^m),$ hence $\xi$ is topologically
isomorphic to the pullback $\eta:=g^{*}\gamma(\mathbf{F}^m)$ (cf. \cite{HuD}, Chapter 3, Theorem 4.7). Since $\xi$
and $\eta$ are topologically isomorphic bundles of class $\mathcal{C}^l,$ then (cf. \cite{Hi}) they are
$\mathcal{C}^l$ isomorphic. Finally, by Proposition \ref{bundle1}, $\eta$ is a $\mathcal{C}^k$ piecewise-algebraic
$\mathbf{F}$-vector
bundle on $L,$ and the proof is complete by Proposition \ref{bundle2}.\qed\vspace*{2mm}\\
{\small\textit{Acknowledgements.} The second named author was partially supported by the National Science Center
(Poland) under Grant number 2018/31/B/ST1/01059.}


{M. Bilski: Faculty of Mathematics and Computer Science, Jagiellonian University, \L ojasiewicza 6, 30-348
Krak\'ow, Poland.\\
E-mail: Marcin.Bilski@im.uj.edu.pl}\vspace*{2mm}\\
{W. Kucharz: Faculty of Mathematics and Computer Science, Jagiellonian University, \L ojasiewicza 6, 30-348
Krak\'ow, Poland.\\
E-mail: Wojciech.Kucharz@im.uj.edu.pl}

\end{document}